\documentclass[letterpaper, 10 pt, conference]{ieeeconf}  % Comment this line out
\IEEEoverridecommandlockouts                              % This command is only
\usepackage{amsmath} % assumes amsmath package installed
\usepackage{amssymb}  % assumes amsmath package installed
\usepackage{graphicx}
\usepackage{subcaption}

\usepackage{algorithm}
\usepackage{algpseudocode}

\newcommand{\vu}{\vec u}
\newcommand{\vz}{\vec z}

\newcommand{\vf}{\vec f}
\newcommand{\vF}{\vec F}
\newcommand{\vY}{\vec Y}
\newcommand{\vU}{\vec U}

\title{\LARGE \bf
Super Resolution for Turbulent Flows in 2D: Stabilized Physics Informed Neural Networks
}

\author{Mykhaylo Zayats, Ma\l gorzata J. Zimo\'{n}, Kyongmin Yeo, Sergiy Zhuk% <-this % stops a space
\thanks{We acknowledge the use of supercomputing facilities provided by the STFC-Hartree Centre.}
\thanks{M. Zayats is with IBM Research - Europe, Dublin, Ireland 
        {\tt\small mykhaylo.zayats1@ibm.com}}%
\thanks{M. J. Zimo\'{n} is with IBM Research - Europe, UK, and The University of Manchester, The School of Mathematics, UK
        {\tt\small malgorzata.zimon@uk.ibm.com}}
\thanks{K. Yeo is with IBM Research
        {\tt\small kyeo@us.ibm.com}}
\thanks{S. Zhuk is with IBM Research - Europe, Dublin, Ireland
        {\tt\small sergiy.zhuk@ie.ibm.com}}
}

\begin{document}

\maketitle
\thispagestyle{empty}
\pagestyle{empty}

%%%%%%%%%%%%%%%%%%%%%%%%%%%%%%%%%%%%%%%%%%%%%%%%%%%%%%%%%%%%%%%%%%%%%%%%%%%%%%%%
\begin{abstract}
We propose a new design of a neural network for solving a zero shot super resolution problem for turbulent flows. We embed Luenberger-type observer into the network's architecture to inform the network of the physics of the process, and to provide error correction and stabilization mechanisms. In addition, to compensate for decrease of observer's performance due to the presence of unknown destabilizing forcing, the network is designed to estimate the contribution of the unknown forcing  implicitly from the data over the course of training. By running a set of numerical experiments, we demonstrate that the proposed network does recover unknown forcing from data and is capable of predicting turbulent flows in high resolution from low resolution noisy observations.
\end{abstract}

%%%%%%%%%%%%%%%%%%%%%%%%%%%%%%%%%%%%%%%%%%%%%%%%%%%%%%%%%%%%%%%%%%%%%%%%%%%%%%%%
\section{INTRODUCTION}
% explain turbulence
Fluid dynamics is a backbone of many contemporary applications in oceanography, weather prediction, energy forecasting, computer graphics, to name a few. Mathematically, an evolution of a fluid flow can be described by the Navier-Stokes equations (NSE), which is a system of nonlinear Partial Differential Equations (PDEs). It is well known that the solution of NSE exhibits complex dynamics even in two dimensions (2D). In particular, the NSE are used as a mathematical model to study turbulence and is related to the concept of deterministic chaos~\cite{Foias15}. Subject to a certain destabilizing forcing, NSE may become very sensitive to small perturbations of initial conditions: initially close-by trajectories diverge over time. In practice, this problem is further amplified since the initial conditions and/or forcing are often not known exactly or can be contaminated with a noise. 

% introducing the problem
For unknown destabilizing forcing, obtaining accurate flow prediction just by solving NSE forward in time is not possible. However, assuming low resolution observations of the flow are available, it may be viable to reconstruct the flow, since observations restrict the set of admissible flows if enough data is provided. Here, we consider the problem of approximating a turbulent flow modelled as a solution of NSE with unknown initial conditions and forcing from low resolution data given in the form of spatial averages over squares $\Omega_j$ covering the computational domain $\Omega$. In the literature, this problem is also referred to as a super resolution problem and more specifically zero-shot super resolution problem, since it does not provide any examples of actual flow snapshots (high resolution ground-truth data)~\cite{Li21}.

\subsection{Observers for Turbulent Systems}
% Observers for turbulent systems
From a control standpoint, a traditional way of estimating states of NSE from observations is to use either a deterministic observer such as Luenberger observer~\cite{Zayats21, Zhuk21, Foias15} or a stochastic filter such as Ensemble Kalman Filter (EnKF)~\cite{Law12}. Both approaches allow for online state estimation, i.e., computing a state estimate once observations at a given time instance become available. In what follows, we will rely upon Luenberger observer, as it has a relatively simple structure (compared to EnKF). In addition, convergence conditions for such observers are available in the literature. Indeed, there have been several attempts to develop sufficient conditions for the convergence of Luenberger observer for NSE. In~\cite{Zayats21}, the authors considered a case of periodic boundary conditions and exact observations taken as spatial averages. These results were improved and extended to a broader class of observation operators, while allowing unknown destabilizing forcing of certain class~\cite{Zhuk21}. The cases of no-slip boundary conditions and periodic boundary conditions together with noisy observations were considered in~\cite{Foias15}. In particular, it was demonstrated that in the presence of bounded noise in observations, the solution of Luenberger observer converges into a certain vicinity of the reference solution.

We stress that in practice, a good guess of the unknown forcing is rarely available, and this can deteriorate practical performance of the observer. To overcome this challenge, we propose a design of a neural network such that the contribution of the unknown forcing is implicitly estimated from the data over the course of training. Hence, for a class of time-independent unknown forcing, our network is expected to perform just as well as if the forcing was known.

\subsection{Physics Informed Neural Networks}
% Explain what is Super-resolution with focus on DNNs
Physics Informed Neural Networks (PINNs) belong to a class of Artificial Neural Networks (ANNs) that explicitly incorporate laws of physics into the network's architecture rather than learning physics purely from data. Development of PINNs for various applications of PDEs is an active research area which has attracted a lot of attention thanks to PINNs ability to extract complex structures from data (see~\cite{Cuomo22} for an extensive overview). PINNs have also been used for solving turbulent NSEs~\cite{Tompson17, WangS20, WangR20} and for solving super-resolution problem for turbulent flows~\cite{Gao21, Jiang20, Li21}.

A typical approach for designing PINN is to use ANN for approximating an unknown function while including PDE as a soft constraint in the loss function that also contains some kind of data misfit~\cite{Cuomo22}. Another strategy is to incorporate parameterized representation of PDE directly into the network structure and use only data misfit as a loss function for optimizing parameters of such a network.

We stress that turbulent flows cannot be exactly predicted even for small errors in forcing/initial conditions~\cite{Foias15,Zhuk21}. This has a severe consequence for designing PINN for turbulent systems. Assume that the time interval is split into training interval $[t_0,t_1]$ and prediction interval $[t_1,t_2]$. By design, PINN informed by the standard NSE acts as a smoother over the training interval and performs well when evaluated on data sampled within the training interval~\cite{Gao21, Li21}. However, due to turbulence, it will likely fail once evaluated using data sampled from the prediction interval. This is demonstrated in~\cite{WangR20} where estimates of turbulent flow are diverging rapidly in prediction interval. Such behaviour is further demonstrated in our own experiments.

In our design, we embed the observer into the network: the dynamic part of the observer's equation (i.e., the part with time derivative) becomes a layer of the network while divergence-free condition is included as a soft constraint in the loss function. The latter  also contains a data misfit term, which is a misfit between low resolution observations and outputs of the network. This new design represents our main contribution: the innovation term of the observer incorporated into the network structure provides an error-correction mechanism based on incoming (in real time) observations. Such mechanism mitigates sensitivity to uncertainty in the initial condition, a key challenge in modelling of turbulent systems, and prevents this uncertainty from being amplified by stabilizing the network. More importantly, the error-correction mechanism is active not only during network training but also in the prediction. This property drastically improves the predictive performance of the designed network for turbulent systems. In what follows, we refer to the designed network as Stabilized PINN or SPINN for short. % it also brings robustness to errors in initial conditions (as is the case for the observer). 

In fact, if the SPINN provides a good numerical approximation for the observer (which is the case if the network is “large enough” as per universal approximation theorem~\cite{Gelenbe99}) the theoretical guarantees of the convergence of the observer (e.g., \cite[Theorem 3.1]{Zhuk21}) apply to the network as well. Our experiments demonstrate that SPINN indeed inherits properties of the Luenberger observer and possesses an error-correction mechanism  which drastically improves its prediction capability compared to PINN informed by the standard NSE.

In addition, our design does not require high resolution data during training, making the approach suitable for a zero-shot super resolution problem. This is in contrast to the network design in~\cite{Jiang20} which requires pairs of low and high resolution snapshots for the network training. 

Finally, we would like to mention an important distinction between SPINNs and conventional numerical methods for solving observers for NSEs: once trained, to compute network's prediction just a simple forward run is required. While the traditional methods require solving the Poisson equation for the pressure gradient computation at every time step which in turn can become quite expensive especially for very fine spatial and temporal discretization.

%\subsection{Key contributions}
%To summarise, our key contributions are the following:
%\begin{enumerate}
%	\item We propose a PINN informed by Luenberger observer for solving a zero-shot %super resolution problem for turbulent flows.
% 	: trained in a self supervised way using only a low resolution data over time interval $[t_0,t_1]$ it is evaluated for predicting high resolution data over time interval $[t_1,t_2]$.
	%\item We demonstrate that the proposed PINN indeed inherits properties of Luenberger observer and possesses an error correction mechanics which drastically improves its prediction capability compared to PINN informed by the standard NSE.
    %\item We demonstrate that Luenberger PINN implicitly recovers unknown forcing from data.
%\end{enumerate}

%%% Local Variables:
%%% mode: latex
%%% TeX-master: "main"
%%% End:

%%%%%%%%%%%%%%%%%%%%%%%%%%%%%%%%%%%%%%%%%%%%%%%%%%%%%%%%%%%%%%%%%%%%%%%%%%%%%%%%
\section{Mathematical Preliminaries}
% \begin{itemize}
%     \item Notation
%     \item NSE 
%     \item Luenberger and convergence conditions (Automatica, Titi for zero dirichlet)
% \end{itemize}
\subsubsection{Notations}
$\mathbb R^n$ denotes the $n$-dimensional Euclidean space. $\Omega$ is a bounded domain in $\mathbb{R}^2$ with boundary $\partial \Omega$.

$L^2(\Omega)$ denotes the space of square-integrable functions on $\Omega$  and $(f,g)=\int_\Omega fgdx$ --  inner product of $L^2(\Omega)$.

\subsubsection{Navier-Stokes equations}
The classical NSE in 2D is a system of nonlinear equations:
\begin{align}
	\label{eq:NSE-system}
	\dfrac{du_1}{dt} - \nu\Delta u_1 + \vu \cdot \nabla u_1 + p_{x} &= f_1,\\
	\dfrac{du_2}{dt} - \nu\Delta u_2 + \vu \cdot \nabla u_2 + p_{y} &= f_2,\\
	\nabla\cdot \vu = (u_1)_{x_1} + (u_2)_{x_2} &= 0,
\end{align}
or in the vector form:
\begin{equation}
	\label{eq:NSE-vector}
	\begin{split}
		\dfrac{d\vu}{dt} -& \nu\Delta u + (\vu \cdot \nabla) \vu + \nabla p = \vf,\\
		&\nabla\cdot\vu = 0.
	\end{split}
\end{equation}
Here, $\vu=(u_1(t,x),u_2(t,x))^\top$ denotes the unknown velocity field and $p(t,x)$ is the unknown scalar pressure field for $(x, t)\in\Omega \times [t_0,\infty)$. $\nu>0$ is a viscosity coefficient and $\vf=[f_1(t,x),f_2(t,x)]^\top$ is a forcing vector. The equation is equipped with the no-slip boundary and initial conditions 
\begin{equation}
    u_i(t,x)=0 \text{ for } x\in\partial\Omega,    
\end{equation}
\begin{equation}
    u_i(0,x)=u_i^0(x),    
\end{equation}
where i=\{1,2\}. In what follows, we assume that the initial conditions and the forcing are such that the Navier-Stokes equations have the unique classical solution. We refer the reader to~\cite{Foias04} for specific existence and uniqueness conditions. 

\subsubsection{Luenberger observer}
One way of estimating $\vz=(z_1(t,x),z_2(t,x))^\top$ of the NSE state $\vu$ from observations in the form 
\begin{equation}
    y_i=Cu_i+\eta_i,
\end{equation}
where $C$ is an observation operator and $\eta_i$ models measurement noise, is to solve the Luenberger observer:
\begin{equation}
	\label{eq:NSE-LO-vector}
	\begin{split}
        \dfrac{d\vz}{dt} - \nu\Delta \vz &+ (\vz \cdot \nabla) \vz + \nabla p = \vec g + \vF, \\
		\nabla\cdot\vz &= 0, \\
		\vz(0,x) &= \vec h(x),
	\end{split}
\end{equation}
where the vector-function $\vec g$ is a guess for the unknown forcing term $\vf$, and $\vec h(x)$ is the guess for the initial condition of NSE, and $\vec F(x,t)$ is the innovation term defined as
\begin{equation}
	\begin{split}
        \vec F(x,t) &= (F_1(x,t),F_2(x,t))^\top, \\
        \quad F_i &=\gamma C^*(y_i - Cz_i).
	\end{split}
\end{equation}

%%% Local Variables:
%%% mode: latex
%%% TeX-master: "main"
%%% End:

%%%%%%%%%%%%%%%%%%%%%%%%%%%%%%%%%%%%%%%%%%%%%%%%%%%%%%%%%%%%%%%%%%%%%%%%%%%%%%%%
\section{PROBLEM STATEMENT}
% \begin{itemize}
% \item assumptions of the noise $\eta_i$??
% \item assumptions on forcing? 
% \item explain that you build DNN architecture by discretizing the observer equation with ZERO forcing, where all the inputs/parameters are known, and you do not discretize  the FORWARD equation with unknowns  (so optimize then discretize rather than the standard discretize then optimize approach)  
% \item then for that discrete equation which preserves the physics you use DNN to estimate grad of p and f by minimizing loss (13), here stress that the network is expected to recover the unknown forcing $f$ from the data 
% \item stress that you enforce a soft version of div-free condition in (13) 
% \item show exact formulas for discretizing laplasian and div z 
% \end{itemize}
Our goal is to design a PINN which solves a \emph{zero-shot super resolution problem for turbulent flows}: given a low resolution and noisy data $\vec Y$ one needs to reconstruct a high resolution approximation $\vU$ of a turbulent fluid flow modelled by means of NSE~\eqref{eq:NSE-vector} with  unknown initial condition and unknown forcing.

In more details, we assume that a 2D viscous fluid flow $\vu=(u_1(t,x),u_2(t,x))^\top$ is modelled as a solution of NSE~\eqref{eq:NSE-vector}, and the latter is subject to unknown initial condition and unknown bounded forcing $\vf$. Moreover, it is assumed that the unknown forcing could be destabilizing and render NSE turbulent for small enough viscosity coefficient $\nu$: i.e., two initially close trajectories diverge over time.

We further assume that the unknown continuous solution $\vu(\vec x,t)$ is related to the data $\vY$ as follows: 
\begin{equation}
	\label{eq:measurements}
\begin{split}
    Y_{i,j}(t)&=\frac{1}{|\Omega_j|}\int_{\Omega_j} u_i(\vec x,t) d\vec x + \eta_i(t)\\ 
    &= (u_i(\cdot,t),b_j) + \eta_i(t), \quad j=1\dots N.
\end{split}
\end{equation}
In other words, $Y_{i,j}(t)$ represents an average of the component of the fluid velocity $u_i$ over non-overlapping subdomains $\Omega_j$ such that $\cup\Omega_j=\Omega$ up to an additive bounded noise $\eta_j$. Here $b_j$ denotes the indicator function of $\Omega_j$ normalized by the measure of $\Omega_j$, $|\Omega_j|$. For example, to give a precise meaning to the low resolution data one can take $\Omega_j$ to be large enough rectangles, covering $\Omega$. In this case, the data $\vY$ is indeed a low resolution piecewise constant (up to noise $\eta_j$) approximation of $\vu$ by means of averages of its components over $\Omega_j$.

The high resolution approximation $\vU$ of $\vu$ refers to approximating the values of the continuous in space vector-function $\vu$ at a given set of grid points $\vec x_1\dots \vec x_K$, densely covering the domain $\Omega$, by $\vU(\vec x_j)$.

%%% Local Variables:
%%% mode: latex
%%% TeX-master: "main"
%%% End:

%%%%%%%%%%%%%%%%%%%%%%%%%%%%%%%%%%%%%%%%%%%%%%%%%%%%%%%%%%%%%%%%%%%%%%%%%%%%%%%%
\section{STABILIZED PINNs}
\label{s:spinn}
To design a PINN which solves the \emph{zero-shot super resolution problem for turbulent flows}, we employ the following optimize-discretize-reconstruct strategy:
\begin{itemize}
    \item [1)] \textit{Optimize:} Luenberger observer~\eqref{eq:NSE-LO-vector} with the innovation term $\vF=(F_1,F_2)^\top $, where
    \begin{equation}\label{eq:innov-term}
        F_i =\gamma\sum_{j=1}^N(Y_{i,j}(t) - (z_i(t,\cdot),b_j)) b_j, 
    \end{equation}
    is used as a meta layer of the network: it provides a mathematical description of the underlying physical process with the error-correction mechanism for treating turbulence. 

    \item [2)] \textit{Discretize:}
    the equations of the observer~\eqref{eq:NSE-LO-vector}, \eqref{eq:innov-term} are then discretized in time by the splitting method and in space by the finite difference method.

    \item [3)] \textit{Reconstruct:} to get high resolution approximation $\vU$, the computation of the sum of pressure gradient and forcing guess $\vec g$ is replaced with ANN parametrization. The constructed PINN is trained by minimizing the corresponding loss function, consisting of data misfit and divergence-free terms.
\end{itemize}

\subsubsection{Optimize}
It was demonstrated in~\cite{Zhuk21} that the Luenberger observer ~\eqref{eq:NSE-LO-vector}, \eqref{eq:innov-term} converges globally for spatial average measurements for the case of known and unknown forcing if the observer gain $\gamma$ and partition $\{\Omega_j\}$ verify conditions presented in~\cite[Theorem 3.1.]{Zhuk21}. These conditions suggest that the convergence is guaranteed once the size of a rectangle $\Omega_j$ is proportional to $\nu^2/\mathrm{log}^\frac{1}{2}\nu$. The conditions in~\cite{Zhuk21}, however, are obtained for periodic boundary conditions and for exact observations. We stress that no-slip boundary conditions considered here would result in even more conservative estimates for $\gamma$ and $\{\Omega_j\}$. The case of observations with bounded deterministic noise, no-slip boundary conditions but known forcing was considered in~\cite{Foias15}. The authors proved convergence of the observer to the true state up to a term which depends on the noise bound. The obtained convergence requirements for $\gamma$ and $\{\Omega_j\}$ in~\cite{Zhuk21} and~\cite{Foias15} are rather conservative. From practical considerations, they serve primarily as a reference point and could be significantly amended, as demonstrated in our experiments.

The ability of the observer to be robust to errors in initial conditions is crucial for modelling a turbulent system and is a much desired property in PINN design. For this reason, in our design of Stabilized PINN we incorporate a discretized version of Luenberger observer into a network structure.

\subsubsection{Discretize} 
We discretize observer~\eqref{eq:NSE-LO-vector} in time by using the splitting method: 
\begin{align}
    \vz^{\,*} &= \vz^{\,t} + \left\{ \nu \Delta \vz^{\,t}  -(\vz^{\,t} \cdot \nabla) \vz^{\,t} + \vec F^{\,t}\right\}\delta t,\label{eqn:NSE-step-1}\\
    \vz^{\,t+1} &= \vz^{\,*} - \delta t \nabla p^{t+1} + \delta t \vec g^{\,t},\label{eqn:NSE-step-2}
\end{align}
in which $\delta t$ is the time step size, the superscript $t$ denotes the time step and $\vz^{\,*}$ is a tentative velocity. Then, the pressure $p^{t+1}$ is obtained by solving the Poisson equation
\begin{equation} \label{eqn:Poisson}
    \Delta p^{t+1} = \frac{1}{\delta t}\nabla \cdot (\vz^{\, *}+\delta t \vec g^{\,t}),
\end{equation}
with the boundary condition,
\begin{equation}
\vec n \cdot \nabla p^{t+1} = 0,~\text{for}~x \in \partial \Omega.
\end{equation}
Here, $\vec n$ denotes an outward normal vector.

For computing $\nabla \cdot \vz^{\,t}$ and $\Delta \vz^{\,t}$  in~\eqref{eqn:NSE-step-1} a spatial discretization is performed on a rectangular uniformly spaced mesh in both $x_1$ and $x_2$ directions: $\delta x_1 = \delta x_2 = \delta x$. The divergence $\nabla \cdot \vz^{\,t}$ and the Laplacian $\Delta \vz^{\,t}$ are approximated by the second order central difference scheme:
\begin{equation}
    \label{eq:fin-diff-approx}
    \frac{\partial z_1^{\,t}}{\partial x_1} \approx \frac{z^{\,t}_{1,i+1}-z^{\,t}_{1,i-1}}{2\delta x},
\end{equation}
\begin{equation}
    \label{eq:div-approx}
    \frac{\partial^2 z_1^{\,t}}{\partial x_1^2} \approx \frac{z^{\,t}_{1,i+1}-2z^{\,t}_{1,i}+z^{\,t}_{1,i-1}}{(\delta x)^2}. 
\end{equation}

\noindent For the nonlinear term $(\vz^{\,t} \cdot \nabla) \vz^{\,t}$ an upwind method is employed, e.g.,
\begin{equation}
    \label{eq:upwind-approx}
    z_1^{\,t}\frac{\partial z_2^{\,t}}{\partial x_1} \approx z^{\,t}_{1,i} \frac{z^{\,t}_{2,i+1}-z^{\,t}_{2,i-1}}{2\delta x} - |z^{\,t}_{1,i}| \frac{z^{\,t}_{2,i+1}-2z^{\,t}_{2,i}+z^{\,t}_{2,i-1}}{2\delta x},
\end{equation}
where the subscript $i$ denotes the spatial grid index, $z_{1,i} = z_1(i\delta x)$.

\subsubsection{Reconstruct}
\begin{figure}
    \centering
    \includegraphics[width=0.48\textwidth]{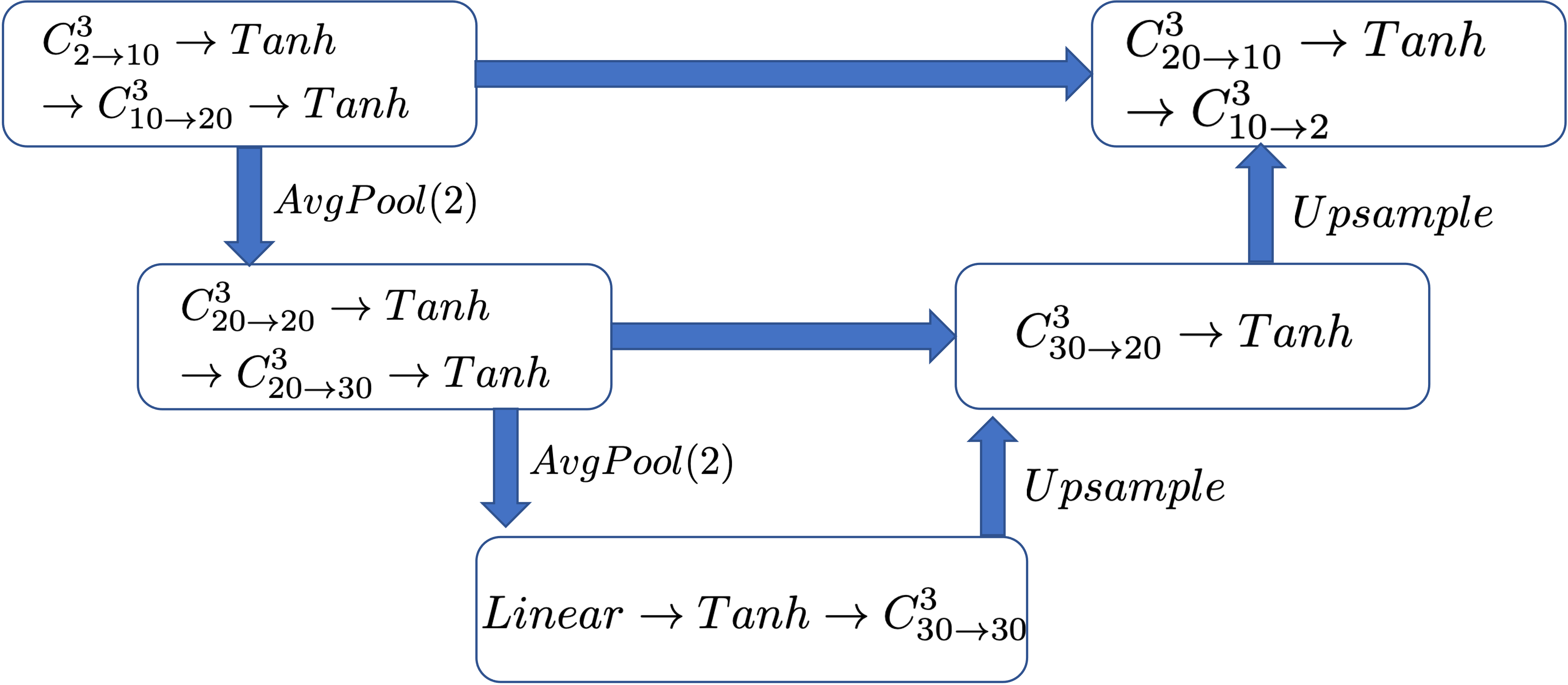}
    \caption{Sketch of the U-Net architecture of $H$. Here, $C^a_{b\rightarrow c}$ denotes a 2D convolution neural network with input channel $b$, output channel $c$, and filter width $a$.}
    \label{fig:unet}
\end{figure}

\begin{algorithm}
    \caption{An algorithm for solving super resolution problem by using SPINN}\label{alg:spinn_alg}
    \begin{algorithmic}
    \Require data $\vY$ over time interval $\left[t_0;T\right]$;
    \State 1. split $\left[t_0;T\right]$ into $\left[t_0;t_1\right]$ and $\left[t_1;T\right]$;
    \State 2. construct SPINN~\eqref{eqn:NSE-step-1}, \eqref{eqn:ANN-step-2};
    \State 3. train SPINN on $\vY_{\left[t_0;t_1\right]}$ by selecting parameters of U-NET by minimizing the loss~\eqref{eq:ann-loss};
    \State 4. predict $U$ by evaluating SPINN on $\vY_{\left[t_1;T\right]}$;
    \State 5. repeat steps 2-4 with different $\gamma$ and select the best SPINN.
\end{algorithmic}
\end{algorithm}

We stress that in practice a good guess $\vec g$ of the unknown forcing $\vf$ is not available, and this constitutes a challenge for observer's performance. Indeed, in the discretization presented above, the error between $\vec g$ and $\vf$ could affect the update of $\vz^{\, t+1}$ from $\vz^{\, t}$ as per~\eqref{eqn:NSE-step-1}, as well as updating pressure as per~\eqref{eqn:Poisson} at every time step. In turn, this can significantly compromise the quality of the high-resolution approximation $\vU$ given by $\vz^{\, t+1}$. Note that this problem is independent of the discretization, as $\vec g$ appears in the continuous equation~\eqref{eq:NSE-LO-vector}. 

Yet another challenge is of computational nature: the need to solve Poisson equation~\eqref{eqn:Poisson} at every time step makes the above scheme quite expensive, especially for very fine spatial and temporal discretization.  

To overcome these two challenges, we employ an ANN $H(\cdot)$ which takes the tentative velocity $\vz^{\,*}$ as an input and computes $-\nabla p^{t+1}+\vec g^{\, t}$ as an output without directly solving~\eqref{eqn:Poisson}. Hence, the sum of pressure gradient and forcing in~\eqref{eqn:NSE-step-2} is replaced by the network $H(\vz^{\,*})$ and results in
\begin{equation} \label{eqn:ANN-step-2}
    \vz^{\,t+1} = \vz^{\,*} + H(\vz^{\,*}) \delta t.   
\end{equation}
There are several suitable network architectures for representing $H(\vz^{\,*})$~\cite{Li21}. Here, we use the U-Net architecture, a well established option for such input-output relationships~\cite{Rnneberger15}. Figure~\ref{fig:unet} shows an outline of the architecture of the U-Net network that we employ. We leave the detailed investigation of the optimal network architecture for future research. 

Equations~\eqref{eqn:NSE-step-1} and~\eqref{eqn:ANN-step-2} define a structure of the SPINN which preserves physics by embedding the Luenberger observer. It is also expected that the SPINN implicitly estimates forcing guess $\vec g$ from data. To train parameters of the SPINN and more specifically weights of the U-Net network $H(\cdot)$, we are minimizing the following objective function:
\begin{align}
    \label{eq:ann-loss}
    \mathcal{L} = \sum_{t=1}^T & \sum_{j} \| \vec Y_{j}(t\delta t) - (\vz^{\,t},b_j) \|^2_2 \nonumber \\
    & + \lambda \| \nabla \cdot \vz^{\,t} \|^2_2.
\end{align}
The first term of $\mathcal{L}$ is the error between a low resolution projection of $\vz^{\,t}$ and the actual low resolution observations $\vY_{j}(t\delta t)$. The second term imposes the soft version of the divergence-free condition on $\vz^{\,t+1}$ and $\lambda$ is a regularization parameter. 

An outline of the approach for solving the super resolution problem for turbulent flows with the proposed SPINN is summarized in Algorithm~\ref{alg:spinn_alg}.

%%% Local Variables:
%%% mode: latex
%%% TeX-master: "main"
%%% End:

%%%%%%%%%%%%%%%%%%%%%%%%%%%%%%%%%%%%%%%%%%%%%%%%%%%%%%%%%%%%%%%%%%%%%%%%%%%%%%%%
\section{EXPERIMENTS}
\begin{figure*}
\centering
  \begin{minipage}[t]{0.45\textwidth}
    \includegraphics[width=\textwidth]{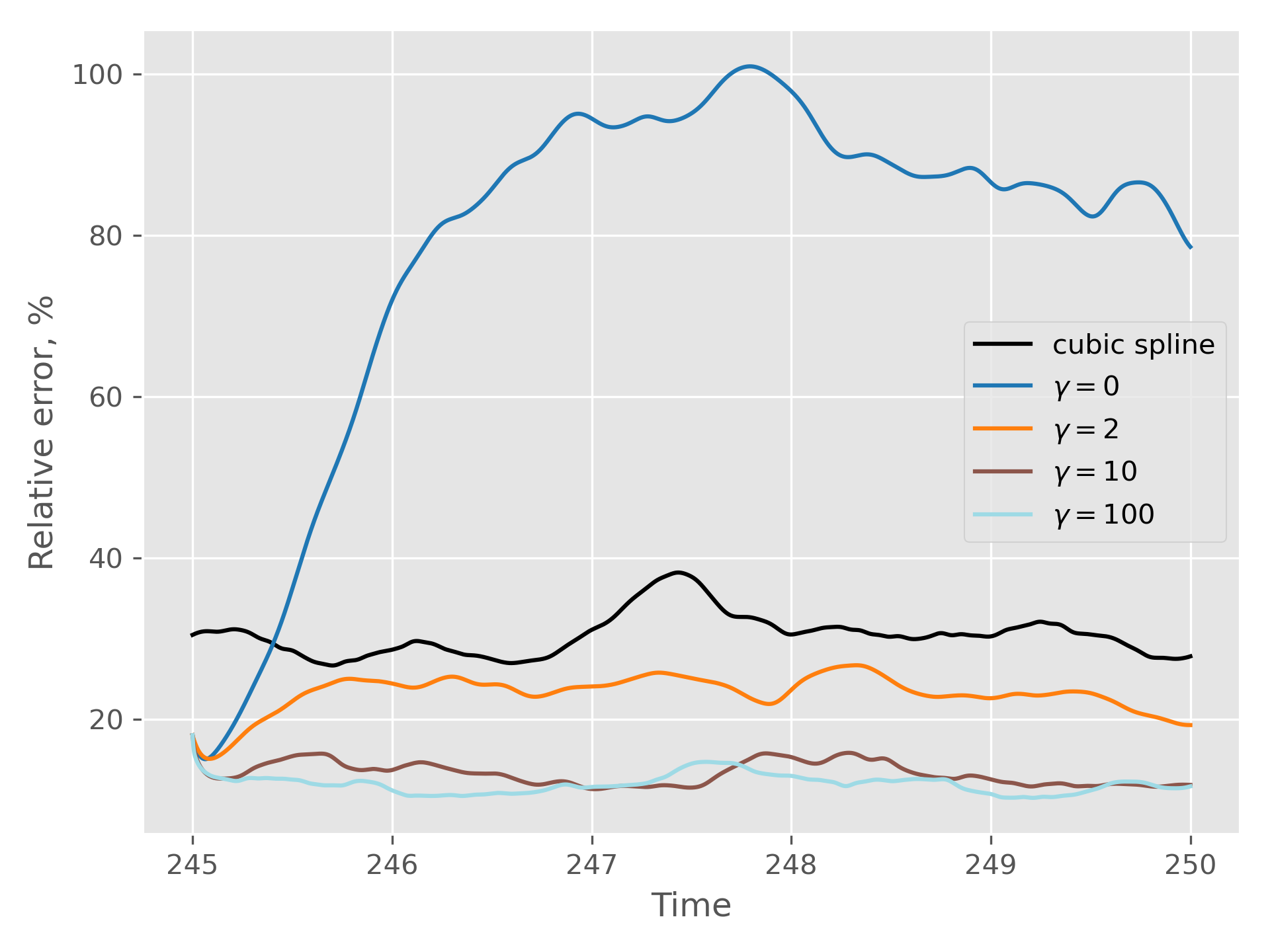}
    \caption{Comparison of prediction error obtained with the SPINN approach ($\gamma=\{2,10,100\}$), PINN approach informed by standard NSE ($\gamma=0$) and bi-cubic interpolation. Computed for the first component of the flow.}
  \label{fig:compare_predictions}
  \end{minipage}
  \quad
  \begin{minipage}[t]{0.45\textwidth}
    \includegraphics[width=\textwidth]{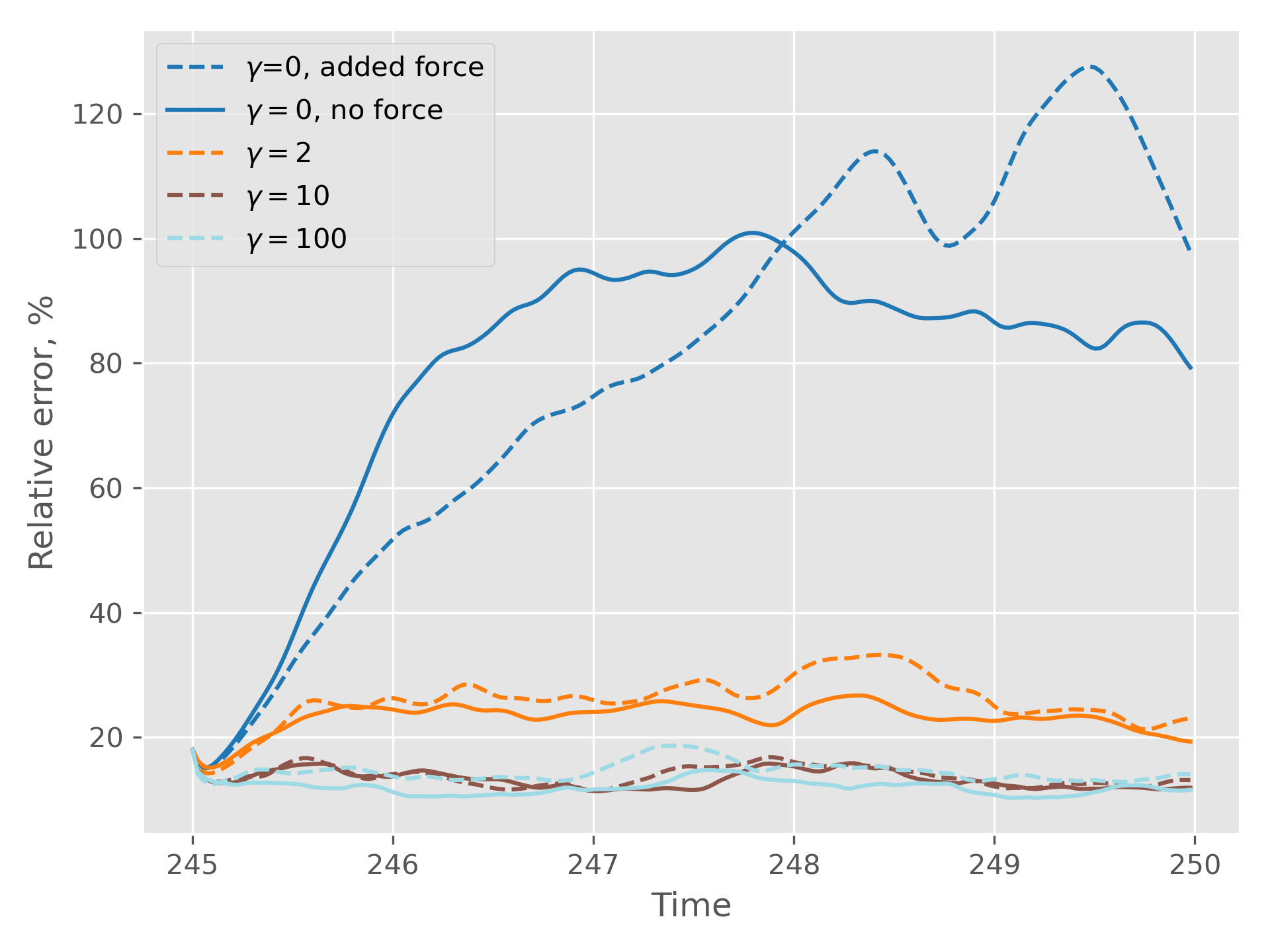}
    \caption{Comparison of prediction error obtained with the SPINN with unknown forcing (solid lines) and known forcing (dashed lines) computed for $\gamma=\{0,2,10,100\}$. Computed for the first component of the flow.}
    \label{fig:comp_forcing}    
  \end{minipage}
\end{figure*}

% \begin{figure}
%   \centering
%     \includegraphics[width=0.49\textwidth]{figures/cr4_errors_no_force_interp.png}
%     \caption{Comparison of prediction error obtained with the SPINN approach ($\gamma=\{2,10,100\}$), PINN approach informed by standard NSE ($\gamma=0$) and bi-cubic interpolation. Computed for the first component of the flow.}
%   \label{fig:compare_predictions}
% \end{figure}

% \begin{figure}
%     \centering
%     \includegraphics[width=0.49\textwidth]{figures/cr4_errors_with_no_force.png}
%     \caption{Comparison of prediction error obtained with the SPINN with unknown forcing (solid lines) and known forcing (dashed lines) computed for $\gamma=\{0,2,10,100\}$. Computed for the first component of the flow.}
%     \label{fig:comp_forcing}
% \end{figure}

\begin{figure*}
  \centering
  \begin{subfigure}[t]{0.45\textwidth}
    \includegraphics[width=\textwidth]{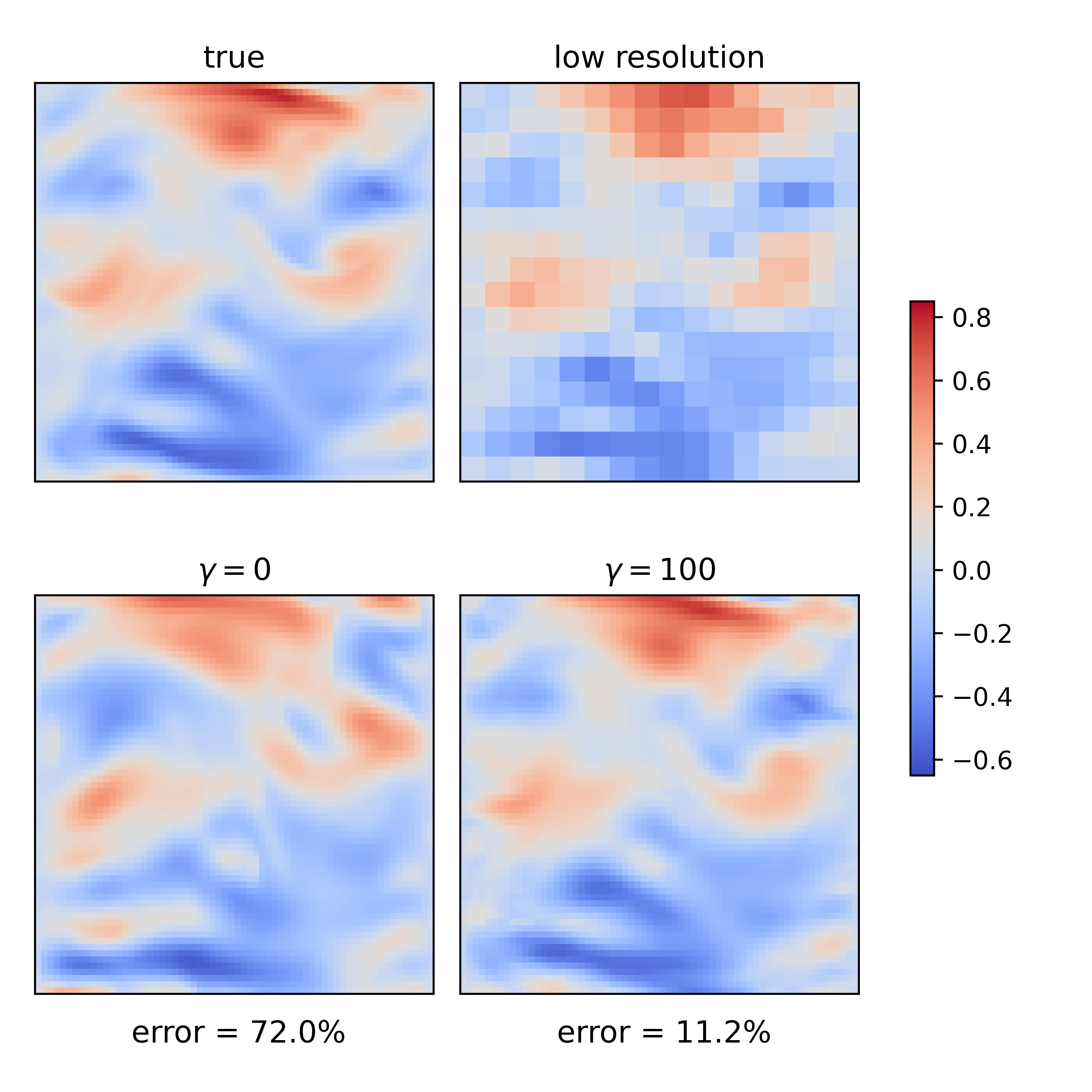}
    \caption{Prediction for $t=246$.}
  \end{subfigure}
  \quad
  \begin{subfigure}[t]{0.45\textwidth}
    \includegraphics[width=\textwidth]{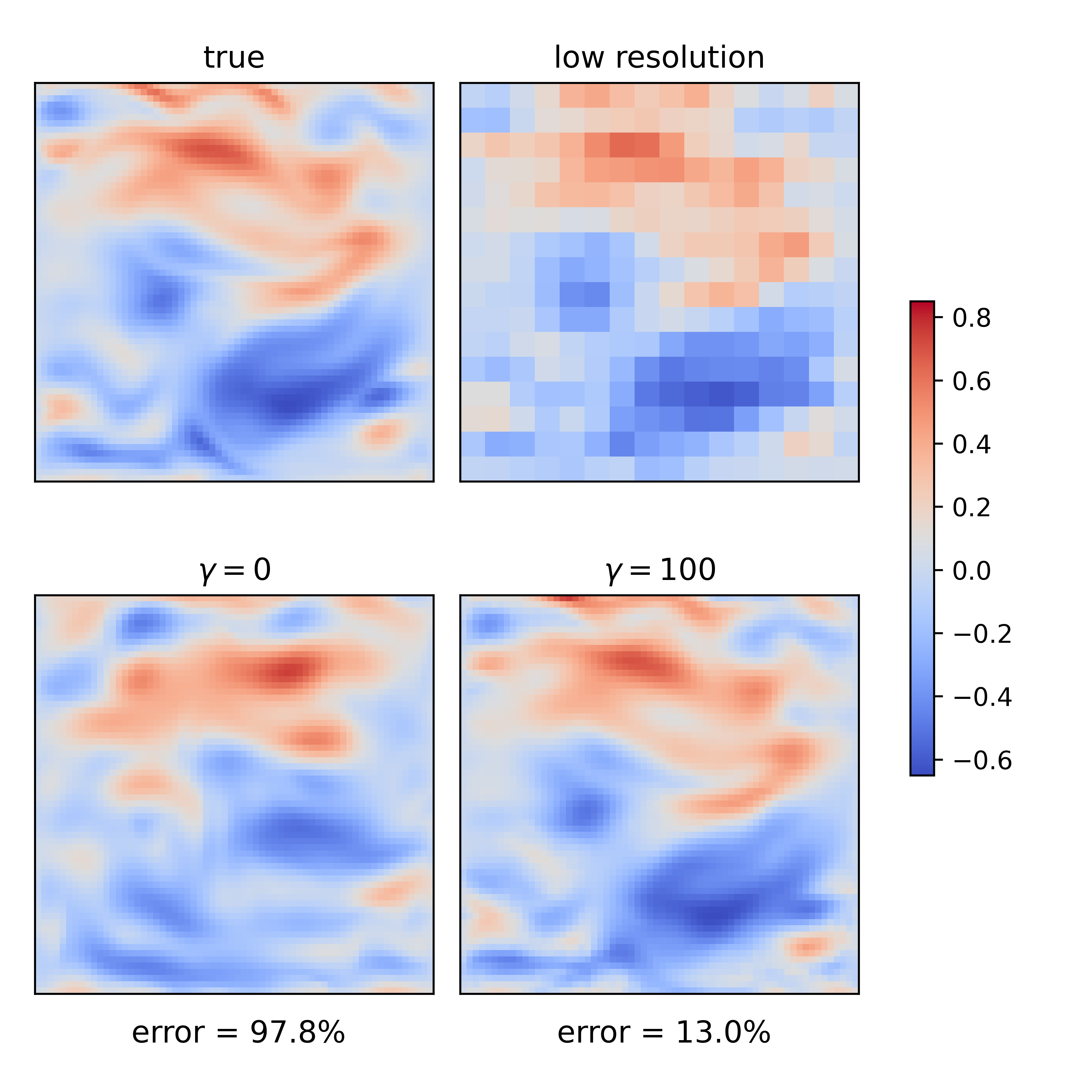}
    \caption{$t=248$.}
  \end{subfigure}
  
  \caption{Examples of the super resolution solutions for $z_1$ at $t=246$ (a) and $t=248$ (b). The plots show the true solution, low resolution input, and super-resolution predictions for $\gamma=0$ and $\gamma=100$.}
  \label{fig:compare_velocity_fields}
\end{figure*}

To demonstrate the efficiency of the proposed approach, we perform a set of numerical experiments. First, we generate observations by running a forward solution of NSE~\eqref{eq:NSE-vector} defined over a spatial domain $\Omega=[0; 2\pi]\times[0; 2\pi]$ and a temporal domain $t=[0; 250]$. The equation is equipped with the viscosity coefficient $\nu=0.01$, the initial conditions taken as an exponential vector function and the forcing taken as 6-th Fourier mode with an amplitude 1. Such configuration is selected intentionally to produce turbulent velocity flow. The reference solution $\vu$ is generated by using finite element method (FEM) implemented using Fenics package~\cite{Alnaes15}. For spatial discretization, we construct a uniform grid with 64 by 64 nodes and generate 7938 linear triangular elements. For temporal discretization, we execute the solver~\cite[section 3.4]{Alnaes15} for 500000 time steps with the time step $\delta t=0.0005$.

FEM solver is computed with $\delta t$, while the measurements are saved every $t_s=5\delta t$, resulting in 100000 snapshots. From high resolution data of the size 64 by 64, we produce low resolution observations through the average pooling with kernel size of (4, 4). We do not add any noise to observations explicitly, however, since those observations are generated from a numerical solution of NSE they contain unknown but bounded FEM approximation errors. The low resolution observations of the velocity field of the size 16 by 16 are then divided into two sets of $N_t=98000$ steps for $t_s \in [0, 245]$ and $N_v=2000$ steps for $t_s \in (245, 250]$ for training and prediction, respectively.

To train the SPINN, we use Adam~\cite{kingma2014adam}, a first-order gradient-based optimization algorithm, with learning rate of $10^{-4}$ and no weights decay. At each stochastic gradient descent step, twenty 200-time-step samples are randomly selected for a mini-batch.

For SPINN training it is also important to select a proper value of $\gamma$. There are theoretical considerations for selecting $\gamma$, however those are not always suitable in practice. For instance, for $\nu=0.01$ and the size of a rectangle $\Omega_j$ corresponding to 16 by 16 rectangles taken as $h^2=(2\pi/16)^2=0.154$ the inequality in~\cite[Theorem 3.1.]{Zhuk21} is infeasible. At the same time, for $\nu=0.01$ that inequality predicts $h^2=0.0155$ which corresponds to 52 by 52 rectangles $\Omega_j$. For this reason, in our experiments we run several SPINNs with different values of $\gamma$ to select the optimal one.

\subsection{Predictive power}

Here, we demonstrate that the proposed SPINN which is informed by the Luenberger observer indeed inherits the property of the observer and, once it is trained, performs the super resolution computation in prediction mode without deprecation of quality over time. To this end, we train the network with different values of the observer gain $\gamma$ and evaluate the models on data sampled from the prediction time interval. Note that the high resolution data is not presented to SPINN in the training.

To assess the quality of network predictions, we introduce the prediction error computed as the relative $l_2$-error between the predicted high resolution solution and the reference solution:
\begin{equation}
 \epsilon_{i,t} = \| u_i(t, .) - z_i(t, .) \|_2 / \|u_i(t, .) \|_2.
\end{equation}

The prediction error is shown in Figure~\ref{fig:compare_predictions}. It demonstrates that for $\gamma=0$, when the innovation term has no impact on the network, i.e., PINN in informed by the standard NSE, the model fails to generalize to the unseen data due to the turbulent behaviour of the underlying process. As a result, the predicted velocity diverges from a reference solution over time. On the other hand, for $\gamma>0$ the innovation term performs error correction from the incoming real time data and updates predicted velocity. This mitigates the impact of turbulence, so that the predictive error for networks trained with $\gamma>0$ converges to reasonably small levels. We note that a good accuracy is already achieved for $\gamma\approx 10$; there is very little improvement for higher values of $\gamma$. This is further confirmed by the results in the Table~\ref{tab:compare_error_for_interp} presenting $\epsilon_{1,t}$ averaged over the time interval $[246,250]$ (errors from time range of [245;246] are excluded to remove the influence of initial phase). Although good results are obtained for $\gamma=100$, we observe that for $\gamma \gg 10$ the innovation term introduces numerical instabilities due to discontinuities of reconstruction of low resolution error resulting in degradation of prediction quality. For small $\gamma$, the impact of the innovation term is not strong enough, resulting in a larger convergence zone compared to a zone produced by optimal $\gamma$.

Examples of the high resolution reconstruction are demonstrated in Figure~\ref{fig:compare_velocity_fields}. It clearly illustrates the ability of the SPINN to reconstruct small-scale features. In contrast, the standard PINN, i.e., $\gamma=0$, deviates from the reference solution due to the lack of the self-correction mechanism.

% comparison against bi-cubic interpolation
To demonstrate that the SPINN indeed incorporates the underlying physics, we compare the SPINN prediction against bi-cubic interpolation of low resolution data. Figure~\ref{fig:compare_predictions} shows how the SPINN with $\gamma>0$ outperforms upscaling of low resolution inputs with spline interpolation. Clearly, without the knowledge about the physics of the process, the interpolation approach leads to over-smoothing of flow features. It is interesting to note, however, that oversimplifying a physical model such as using NSE with perturbed inputs for turbulent flows produces results that are less accurate even compared to the physics unaware interpolation method.

\subsection{Forcing reconstruction}
To demonstrate the ability of the SPINN to recover unknown forcing $\vf$ we compare it against SPINN for which the estimate of forcing is known exactly, i.e., $\vec g=\vf$ and $\vf$ is defined as 6-th Fourier mode. For this, we modify equation~\eqref{eqn:ANN-step-2} in the following way
\begin{equation} \label{eqn:ANN-step-2-forcing}
    \vz^{\,t+1} = \vz^{\,*} + H(\vz^{\,*}) \delta t + \vf \delta t. 
\end{equation}

Fig.~\ref{fig:comp_forcing} shows $\epsilon_{1,t}$ obtained with the Luenberger PINNs with unknown and explicitly defined forcing. It is found that, regardless of $\gamma$, the errors with and without the explicit forcing term are similar to each other, confirming the ability of the network to implicitly recover forcing $\vf$ from data. For a quantitative comparison, the averaged $\epsilon_{1,t}$ is provided in Table~\ref{tab:compare_error_for_interp}. We also note that for 3 cases of $\gamma$ presented in the table, the network with unknown forcing slightly outperforms the network with known forcing. This is likely due to the complex interplay between approximation errors introduced by discretization and network optimisation.
% This is likely due to the impact of approximation errors.

\begin{table}[t]
\begin{tabular}{l|llll}
  forcing  & $\gamma=2$  & $\gamma=10$ & $\gamma=100$ & cubic spline \\  \hline
  unknown  & 23.74 $\%$  & 13.04 $\%$  & 11.90 $\%$    & 30.86 $\%$ \\
  known    & 26.99 $\%$  & 13.60 $\%$  & 14.56 $\%$        & -    
\end{tabular}
\caption{Average relative errors obtained for Luenberger PINN with unknown and known forcing and different values of $\gamma$ and bi-cubic interpolation.}
\label{tab:compare_error_for_interp}
\end{table}

\begin{table}
\centering
\begin{tabular}{l|llll}
  & $N_t=98000$ &  $75\% N_t$ & $50\% N_t$  & $25\% N_t$  \\
  \hline
  $\gamma=10$ &  13.04 $\%$&  14.23 $\%$&  15.95$\%$ & 20.10 $\%$\\
  $\gamma=100$ &  11.90 $\%$&  11.12 $\%$&  13.68$\%$ & 14.31 $\%$ 
\end{tabular}
\caption{Averaged predictive errors of SPINN trained on datasets generated using the last $75\% N_t$, $50\% N_t$, $25\% N_t$ of samples from the original dataset $N_t=98000$ and evaluated on the same prediction dataset.}
\label{tab:compare_error_for_diff_Nt}
\end{table}

\begin{table}
\centering
\begin{tabular}{l|lll}
               & $N_t=98000$ & $50\% N_t$ & $25\% N_t$  \\ \hline
  $\gamma=10$  & 13.04 $\%$  & 19.07 $\%$ & 23.65$\%$\\
  $\gamma=100$ & 11.90 $\%$  & 12.84$\%$  & 13.96 $\%$
\end{tabular}
\caption{Averaged prediction errors of SPINN trained on datasets generated by taking every 2nd and 4th sample from the original dataset and evaluated on the same prediction dataset.}
\label{tab:compare_error_for_diff_Nt2}
\end{table}

\subsection{Generalization on smaller dataset}

When dealing with ANNs, a significant challenge is the necessity of a large amount of data for training. In this test, we analyse the impact of the amount of training data on the performance of the proposed SPINN. We consider two scenarios: in the first, we reduce the number of training samples $N_t$ by taking the last $75\%$, $50\%$ and $25\%$ of samples, resulting in data size of $N_t=73500$,  $N_t=49000$,  $N_t=24500$, respectively. In the second scenario, the time interval between two consecutive snapshots is increased to $2t_s$ and $4t_s$, resulting in $N_t=49000$,  $N_t=24500$ samples, respectively. For each case, we train the network on a smaller dataset and the errors are evaluated using the $N_v=2000$ samples. The average errors of approximations are reported in Table~\ref{tab:compare_error_for_diff_Nt} and Table~\ref{tab:compare_error_for_diff_Nt2} for the first and second scenario, respectively. Although both experiments demonstrate an increase of averaged prediction error, for the SPINN model with $\gamma=100$ such an increase is only marginal: from $11.9\%$ to $14.31\%$ for the first scenario and to $13.96\%$ for the second. For the SPINN model with $\gamma=10$ an increase is more significant: from $13.4\%$ to $20.1\%$ for the first scenario and from $13.4\%$ to $23.65\%$ for the second. We note that in any case averaged prediction error of the SPINN models is lower than of the bi-cubic interpolation which is $30.86\%$ confirming the ability of the PINN to generalize on a smaller training data. 

The difference in the patterns of averaged prediction error growth between models with $\gamma=100$ and $\gamma=10$ demonstrates the impact of the innovation term on network training, in particular for training with less data. Training on smaller datasets results in less error correction during the training which, as the results in the Table~\ref{tab:compare_error_for_diff_Nt} and Table~\ref{tab:compare_error_for_diff_Nt2} suggest, may be compensated by larger values of $\gamma$.

%%% Local Variables:
%%% mode: latex
%%% TeX-master: "main"
%%% End:

%%%%%%%%%%%%%%%%%%%%%%%%%%%%%%%%%%%%%%%%%%%%%%%%%%%%%%%%%%%%%%%%%%%%%%%%%%%%%%%%
\section{CONCLUSIONS}
A stabilized PINN was proposed for solving a zero-shot super resolution problem for turbulent flows. It was demonstrated experimentally that: 
\begin{itemize}
    \item the proposed network indeed inherits properties of Luenberger observer and possesses an error-correction mechanism which drastically improves its prediction capability compared to PINN informed by the standard NSE,
    \item the network implicitly recovers unknown forcing from data.
\end{itemize}

The designed network has the potential of integration of parameter estimation (e.g., unknown viscosity) without a drastic change of the network's architecture -- an interesting subject for future research. 

%%%%%%%%%%%%%%%%%%%%%%%%%%%%%%%%%%%%%%%%%%%%%%%%%%%%%%%%%%%%%%%%%%%%%%%%%%%%%%%%

\end{document}